\documentclass[a4paper,11pt]{article}
\usepackage[latin1]{inputenc}
\usepackage[francais,english]{babel}    
\usepackage{amssymb}
\usepackage{amsmath}
\usepackage{amsthm}
\usepackage{multicol}
\usepackage{graphicx}
\usepackage{color}
\usepackage[T1]{fontenc}
\usepackage{yfonts}
\usepackage{placeins}
\usepackage{array}
\usepackage{dsfont}
\usepackage{stmaryrd}
\usepackage{verbatim}
\usepackage{caption}
\usepackage{transparent}
\usepackage{bbold}
\usepackage{hyperref}
\usepackage{enumitem}
\usepackage{tikz}
\usetikzlibrary{patterns}
\usepackage{tikz-3dplot}

\newtheorem{thm}{Theorem}[section]
\newcommand{\bthm}{\begin{thm}}
\newcommand{\ethm}{\end{thm}}

\newtheorem{thmi}{Theorem}
\newcommand{\bthmi}{\begin{thmi}}
\newcommand{\ethmi}{\end{thmi}}

\newtheorem{cori}[thmi]{Corollary}
\newcommand{\bcori}{\begin{cori}}
\newcommand{\ecori}{\end{cori}}

\newtheorem{mthm}{Theorem}
\newcommand{\bmthm}{\begin{mthm}}
\newcommand{\emthm}{\end{mthm}}

\newtheorem{mcor}[mthm]{Corollary}
\newcommand{\bmcor}{\begin{mcor}}
\newcommand{\emcor}{\end{mcor}}
\newtheorem{mconj}[mthm]{Conjecture}
\newcommand{\bmconj}{\begin{mconj}}
\newcommand{\emconj}{\end{mconj}}
\newtheorem{mpro}[mthm]{Proposition}
\newcommand{\bmpro}{\begin{mpro}}
\newcommand{\empro}{\end{mpro}}

\newtheorem*{conj}{Conjecture}
\newcommand{\bconj}{\begin{conj}}
\newcommand{\econj}{\end{conj}}

\newtheorem*{question}{Question}
\newcommand{\bq}{\begin{question}}
\newcommand{\eq}{\end{question}}

\newtheorem*{thn}{Theorem}
\newcommand{\bthn}{\begin{thn}}
\newcommand{\ethn}{\end{thn}}

\newtheorem{exo}{Exercise}
\newcommand{\bex}{\begin{exo}}
\newcommand{\eex}{\end{exo}}

\newtheorem{sol}{Solution}
\newcommand{\bsol}{\begin{sol}}
\newcommand{\esol}{\end{sol}}

\newtheorem{pro}[thm]{Proposition}
\newcommand{\bpro}{\begin{pro}}
\newcommand{\epro}{\end{pro}}

\newtheorem{cor}[thm]{Corollary}
\newcommand{\bcor}{\begin{cor}}
\newcommand{\ecor}{\end{cor}}

\newtheorem{lem}[thm]{Lemma}
\newcommand{\blem}{\begin{lem}}
\newcommand{\elem}{\end{lem}}

\theoremstyle{definition}

\newtheorem{defi}[thm]{Definition}
\newcommand{\bdf}{\begin{defi}}
\newcommand{\edf}{\end{defi}}

\newtheorem*{defis}{Definition}
\newcommand{\bdfs}{\begin{defis}}
\newcommand{\edfs}{\end{defis}}

\newtheorem*{rmk}{Remark}
\newcommand{\brk}{\begin{rmk} \upshape}
\newcommand{\erk}{\end{rmk}}

\newtheorem*{rmks}{Remarks}
\newcommand{\brks}{\begin{rmks} \upshape}
\newcommand{\erks}{\end{rmks}}

\newtheorem*{exe}{Example}
\newcommand{\bexe}{\begin{exe} \upshape}
\newcommand{\eexe}{\end{exe}}

\newtheorem*{exes}{Examples}
\newcommand{\bexes}{\begin{exes} \upshape}
\newcommand{\eexes}{\end{exes}}

\newtheorem*{pre}{Proof}
\newcommand{\bp}{\begin{pre} \upshape}
\newcommand{\ep}{\hfill \qed \end{pre}}
\newcommand{\epp}{\end{pre}}

\newcommand{\beq}{\begin{eqnarray*}}
\newcommand{\eeq}{\end{eqnarray*}}

\newcommand{\beqn}{\begin{equation}}
\newcommand{\eeqn}{\end{equation}}

\newcommand{\ben}{\begin{enumerate}}
\newcommand{\een}{\end{enumerate}}
\newcommand{\bit}{\begin{itemize} \renewcommand{\labelitemi}{$\bullet$} \renewcommand{\labelitemii}{$\star$}}
\newcommand{\eit}{\end{itemize}}

\newcommand{\ds}{\displaystyle}

\newcommand{\bfg}{
\begin{figure}[H]
\begin{center}}
\newcommand{\efg}{
\end{center}
\end{figure}
\FloatBarrier}

\newcolumntype{M}[1]{>{\raggedright}m{#1}}

\newcommand{\R}{\mathbb{R}}

\newcommand{\N}{\mathbb{N}}
\newcommand{\Z}{\mathbb{Z}}
\newcommand{\C}{\mathbb{C}}

\renewcommand{\O}{\mathbb{O}}

\newcommand{\K}{\mathbb{K}}

\renewcommand{\H}{\mathbb{H}}
\newcommand{\U}{\operatorname{U}}

\newcommand{\bs}{\symbol{92}}

\newcommand{\ov}{\overline}

\renewcommand{\tilde}{\widetilde}

\renewcommand{\dim}{\operatorname{dim}}

\newcommand{\Sp}{\operatorname{Sp}}

\newcommand{\eps}{\varepsilon}
\newcommand{\st}{\, | \,}

\newcommand{\ra}{\rightarrow}

\newcommand{\liml}{\lim\limits}

\newcommand{\f}{\frac}

\renewcommand{\geq}{\geqslant}
\renewcommand{\leq}{\leqslant}

\renewcommand{\log}{\operatorname{log}}

\newcommand{\GL}{\operatorname{GL}}
\newcommand{\SL}{\operatorname{SL}}
\newcommand{\SO}{\operatorname{SO}}
\newcommand{\SU}{\operatorname{SU}}

\newcommand{\mk}{\medskip}

\newcommand{\sign}{\begin{flushright}
Thomas Haettel \\
IMAG, Univ Montpellier, CNRS, France \\
thomas.haettel@umontpellier.fr
\end{flushright}}

\makeatletter
\def\Ddots{\mathinner{\mkern1mu\raise\p@
\vbox{\kern7\p@\hbox{.}}\mkern2mu
\raise4\p@\hbox{.}\mkern2mu\raise7\p@\hbox{.}\mkern1mu}}
\makeatother

\makeatletter
\def\maketitles{%
  \null
  \thispagestyle{empty}%
  \vfill
  \begin{center}\leavevmode
    \normalfont
    {\LARGE \@title\par}%
    \vskip 1.2cm
    {\large \@author\par}%
    \vskip 1.2cm
    {\large \@subtitle\par}%
    \vskip 0.8cm
    {\large \@date\par}%
  \end{center}%
  \vfill
  \null
  \cleardoublepage
  }

\def\date#1{\def\@date{#1}}
\def\author#1{\def\@author{#1}}
\def\title#1{\def\@title{#1}}
\def\subtitle#1{\def\@subtitle{#1}}
\makeatother

\usetikzlibrary{arrows}

\title{Injective metrics on buildings and symmetric spaces}
\author{Thomas Haettel}
\date{\today}

\begin{document}

\selectlanguage{english}

\maketitle

\begin{center}
\begin{minipage}{0.8\textwidth}
\textsc{Abstract.} In this article, we show that the Goldman-Iwahori metric on the space of all norms on a fixed vector space satisfies the Helly property for balls.

On the non-Archimedean side, we deduce that most classical Bruhat-Tits buildings may be endowed with a natural piecewise $\ell^\infty$ metric which is injective. We also prove that most classical semisimple groups over non-Archimedean local fields act properly and cocompactly on Helly graphs. This gives another proof of biautomaticity for their uniform lattices.

On the Archimedean side, we deduce that most classical symmetric spaces of non-compact type may be endowed with a natural invariant Finsler metric, restricting to an $\ell^\infty$ metric on each flat, which is coarsely injective. We also prove that most classical semisimple groups over Archimedean local fields act properly and cocompactly on injective metric spaces. We identify the injective hull of the symmetric space of $\GL(n,\R)$ as the space of all norms on $\R^n$.

The only exception is the special linear group: if $n=3$ or $n \geq 5$ and $\K$ is a local field, we show that $\SL(n,\K)$ does not act properly and coboundedly on an injective metric space.
\end{minipage}
\end{center}

\let\thefootnote\relax

\section*{Introduction}

In this article, we are interested in the relationship between symmetric spaces of non-compact type and Euclidean buildings, on one side, and injective metric spaces and Helly graphs, on the other side.

\mk

A geodesic metric space is called injective if the family of closed balls satisfies the Helly property, i.e. any family of pairwise intersecting balls has a non-empty global intersection. An injective metric space satisfies some properties of nonpositive curvature: it is contractible, any finite group action has a fixed point, and it has a conical geodesic bicombing. One key feature of injective metric spaces is that any metric space embeds isometrically in an essentially unique smallest injective metric space, called the injective hull. Injective metric spaces in geometric group theory have been notably popularized by Lang, who proved that any Gromov-hyperbolic group acts properly and cocompactly on an injective metric space, the injective hull of a Cayley graph (see~\cite[Theorem~1.4]{lang}).

\mk

A geodesic metric space is called coarsely injective if any family of pairwise intersecting balls has a non-empty global intersection, up to increasing the radii by a uniform amount. If a finitely generated group acts properly and cocompactly on a coarsely injective metric space, we can deduce that it is semi-hyperbolic in the sense of Alonso-Bridson. This strategy has been used by Hoda, Petyt and the author to prove that any hierarchically hyperbolic group, including any mapping class group of a surface, is coarsely injective and semi-hyperbolic.

\mk

The discrete analogue of injective metric spaces is the notion of Helly graphs: a connected graph is called Helly if the family of combinatorial balls satisfies the Helly property. The reader is referred to~\cite{helly_groups} for the study of group actions on Helly graphs. One notable result is that a discrete group acting properly and cocompactly on a locally finite Helly graph is biautomatic (see~\cite[Theorem~1.5]{helly_groups}).

\mk

Symmetric spaces of non-compact type and Euclidean buildings already have a CAT(0) metric. Nevertheless, looking for injective metrics on those spaces may provide extra structure. For instance, deciding which CAT(0) groups are biautomatic is very subtle, as Leary and Minasyan recently provided the first counter-examples (see~\cite{leary_minasyan}). On the other hand, any Helly group is biautomatic.

\mk

Our work is based on a very simple remark that, given any set of norms on a vector space satisfying simple conditions, the Goldman-Iwahori metric satisfies the Helly property for closed balls (see~\cite{goldman_iwahori}). The fact that the metric is geodesic will be verified in concrete examples.

\bmpro[Proposition~\ref{pro:helly_norms}] \label{pro:helly_norms_intro}
Let $\K$ denote a valued field, let $V$ denote a $\K$-vector space, and let $X$ denote a set of norms on $V$ satisfying simple conditions (see Proposition~\ref{pro:helly_norms}). For any two elements $\eta,\eta'$ in $X$, let us define the Goldman-Iwahori metric
$$d(\eta,\eta') = \sup_{v \in V \bs \{0\}} \left| \log \f{\eta(v)}{\eta'(v)} \right|.$$
The family of closed balls in the metric space $(X,d)$ satisfies the Helly property.
\empro

\subsection*{Bruhat-Tits buildings}

The first example to which Proposition~\ref{pro:helly_norms_intro} applies is the Goldman-Iwahori space of all ultrametrics norms (see~\cite{goldman_iwahori}). It identifies with the Bruhat-Tits extended building of $\GL(n,\K)$, where $\K$ is a non-Archimedean valued field which is locally compact, or more generally spherically complete. Recall that the Bruhat-Tits building $\ov{X}$ of $\SL(n,\K)$ can be described as the set of all homothety classes of ultrametric norms on $\K^n$ (see~\cite{parreau_immeubles} for instance), and the Bruhat-Tits extended building $X$ of $\GL(n,\K)$ can be described as the set of all ultrametric norms on $\K^n$, also called the Goldman-Iwahori space. Each apartment in $X$ naturally identifies with $\R^n$, and the Goldman-Iwahori metric from Proposition~\ref{pro:helly_norms_intro} is the length metric associated to the standard piecewise $\ell^\infty$ metric on each apartment. We therefore have the following.

\bmthm[Theorem~\ref{thm:BTbuilding_injective}] \label{thm:gln_padic}
Let $\K$ denote any non-Archimedean valued field $\K$ which is spherically complete, and consider the extended Bruhat-Tits building $X$ of $\GL(n,\K)$. Endow $X$ with the Goldman-Iwahori metric, i.e. the length metric associated to the standard piecewise $\ell^\infty$ metric on each apartment. Then $(X,d)$ is injective.
\emthm

Note that a particular case of this result, when the valuation is discrete and the building is simplicial, was already known, combining works of Hirai and Chalopin et al.

\bthn[\cite{hirai_uniform_modular},\cite{chalopin_chepoi_hirai_osajda}]
Let $X$ denote any extended Euclidean building of type $\tilde{A}_{n-1}$. Endow $X$ with the length metric associated to the standard piecewise $\ell^\infty$ metric on each apartment. Then $(X,d)$ is injective.  
\ethn

Our work has the advantage of being valid for a possibly non-discrete valuation if the field $\K$ is spherically complete, and furthermore our proof is extremely simple.

\mk

We can also wonder whether we can apply it to find a Helly graph related to Euclidean buildings. This is indeed the case.

\bmthm[Theorem~\ref{thm:gln_helly}]
Let $\K$ denote any non-Archimedean discretely valued field $\K$, and consider the extended Bruhat-Tits building $X$ of $\GL(n,\K)$. Then the thickening of the vertex set $X^{(0)}$ of $X$ is a Helly graph. In particular, $\GL(n,\K)$ acts properly and cocompactly by automorphisms on a Helly graph.
\emthm

The thickening of $X^{(0)}$ is the graph with vertex set $X^{(0)}$, and with an edge between two vertices if they are at $\ell^\infty$ distance $1$ in some apartment.

\mk

For other classical groups, we can in fact deduce similar results using an embedding in $\GL(n,\K)$.

\bmcor[Theorems~\ref{thm:classical_injective} and \ref{thm:classical_helly}]
Let $\K$ denote a local field of characteristic different from $2$, and let $G$ denote a classical connected semisimple group over $\K$, realized as the identity component of the fixed point set of an involution in the general linear group $\GL(n,\K)$. Then the Bruhat-Tits building of $G$, endowed with the length metric induced from the $\ell^\infty$ metric on the extended Bruhat-Tits building of $\GL(n,\K)$, is injective. Furthermore, the group $G$ acts properly and cocompactly by automorphisms on a locally finite Helly graph.
\emcor

Note that Chalopin et al. proved that any cocompact lattice in a Euclidean building of type $\tilde{C}_n$ acts properly and cocompactly on a Helly graph (see~\cite[Corollary~6.2]{helly_groups}).

\mk

We also easily deduce a result for all classical semisimple Lie groups and their cocompact lattices.

\bmcor[Corollary~\ref{cor:classical_padic}] \label{cor:classical_padic_intro}
Let $G$ denote a classical reductive Lie group over a non-Archimedean local field of characteristic different from $2$, and let $a \geq 0$ denote the number of semisimple factors of type $A$. Then $G \times \Z^a$ acts properly and cocompactly by automorphisms on a locally finite Helly graph.

For any cocompact lattice $\Gamma$ in $G$, the group $\Gamma \times \Z^a$ acts properly and cocompactly by automorphisms on a locally finite Helly graph, and the group $\Gamma$ is biautomatic.
\emcor

Note that Swiatkowski proved that any group acting properly and cocompactly on any Euclidean building is biautomatic (see~\cite[Theorem~6.1]{swiatkowski_biautomatic}). Nevertheless, this provides another perspective on this result.

\subsection*{Symmetric spaces}

The second example to which Proposition~\ref{pro:helly_norms_intro} applies is the symmetric space $X = \GL(n,\R)/\O(n)$ of $\GL(n,\R)$, which may be described as the space of all Euclidean norms on $\R^n$. However, it does not apply directly, since the supremum of two Euclidean norms is no longer Euclidean. So we rather consider the space $\hat{X}$ of all norms on $\R^n$, and use the John-L\"owner ellipsoid to show that $X$ is cobounded in $\hat{X}$.

\bmthm[Theorem~\ref{thm:real_coarse_helly}] \label{thm:gln_real}
Let $X=\GL(n,\R)/\O(n)$ denote the symmetric space of $\GL(n,\R)$, and endow $X$ with the Finsler length metric associated to the standard $\ell^\infty$ metric on each apartment. The injective hull of $X$ is the space $\hat{X}$ of all norms on $\R^n$. Moreover, $X$ is cobounded in $\hat{X}$, which is is proper. As a consequence, $\GL(n,\K)$ acts properly and cocompactly on the injective space $\hat{X}$.\emthm

For other classical groups, we can in fact deduce similar results using an embedding in $\GL(n,\R)$.

\bmthm[Theorem~\ref{thm:classical_coarse_helly}]
Let $G$ denote a classical semisimple non-compact real Lie group which is not of type $\SL$, and let $X$ denote its symmetric space. Then $X$ has a natural Finsler length metric $d$ such that $(X,d)$ is coarsely injective, and its injective hull is proper. In particular, $G$ acts properly and cocompactly by isometries on an injective metric space.
\emthm

We also easily deduce a result for all classical semisimple Lie groups and their cocompact lattices.

\bmcor[Corollary~\ref{cor:classical_real_reductive}] \label{cor:classical_real}
Let $G$ denote any reductive real Lie group, with classical non-compact semisimple factors. Let $a \geq 0$ denote the number of semisimple factors of type $\SL$. Then $G \times \R^a$ acts properly and cocompactly on an injective metric space. In particular, for any cocompact lattice $\Gamma$ in $G$, the group $\Gamma \times \Z^a$ acts properly and cocompactly on an injective metric space. \emcor

Recall that Chalopin et al. proved that any Helly group is biautomatic. This motivates the question whether the non-discrete analogue of this result holds:

\bq
Assume that a finitely generated group $\Gamma$ acts properly and cocompactly on an injective metric space. Is $\Gamma$ biautomatic ?
\eq

\subsection*{The special linear group}

We now turn to the special linear group. According to Theorems~\ref{thm:gln_padic} and \ref{thm:gln_real}, if $\K$ is a local field, we have seen that $\GL(n,\K)$ acts properly and cocompactly on an injective metric space. It is natural to ask what happens for $\SL(n,\K)$. Inspired by the work of Hoda on crystallographic Helly groups (see~\cite{hoda:crystallographic}), we prove the following.

\bmthm[Theorem~\ref{thm:sln_not_helly}]
Let $\K$ be a local field (with characteric different from $2$ if $\K$ is non-Archimedean), and let $n = 3$ or $n \geq 5$. Then $\SL(n,\K)$ is not coarsely injective: $\SL(n,\K)$ does not act properly and coboundedly on an injective metric space.
\emthm

This is also evidence that cocompact lattices in $\SL(n,\K)$ are not expected to be coarsely injective.

\mk

\subsection*{Structure of the article} 

In Section~\ref{sec:injective}, we review the notions of injective metric spaces, Helly graphs and group actions. In Section~\ref{sec:injective_distance}, we present Proposition~\ref{pro:helly_norms} stating that the Goldman-Iwahori metric on the space of all norms satisfies a Helly property for balls. In Section~\ref{sec:buildings}, we apply this construction to Bruhat-Tits buildings, and in Section~\ref{sec:symmetric_spaces}, we apply it to symmetric spaces of non-compact type. In the final Section~\ref{sec:sln_not_helly}, we prove that the special linear group is not coarsely injective.

\mk

\textbf{Acknowledgments:} We would like to thank Victor Chepoi, Bruno Duchesne, Fran\c{c}ois Fillastre, Elia Fioravanti, Anthony Genevois, Hiroshi Hirai, Nima Hoda, Vladimir Kovalchuk, Linus Kramer, Urs Lang, Damian Osajda, Harry Petyt, Betrand R\'{e}my and Constantin Vernicos for interesting discussions and remarks on the first version of the article. We also thank the anonymous referee for interesting comments that helped improve and correct the presentation.

\section{Injective metric spaces and Helly graphs} \label{sec:injective}

In this section, we recall some basic definitions about injective metric spaces and Helly graphs. We refer the reader to~\cite{lang} and \cite{helly_groups} for more details.

\mk

A metric space $(X,d)$ is called \emph{injective} if, for any family $(x_i)_{i \in I}$ of points in $X$ and $(r_i)_{i \in \N}$ of nonnegative real numbers satisfying
$$ \forall i,j \in I, r_i+r_j \geq d(x_i,x_j),$$
the family of balls $(B(x_i,r_i))_{i \in \N}$ has a non-empty global intersection.

\mk

In case the metric space $(X,d)$ is geodesic, it is injective of and only if the family of balls satisfy the \emph{Helly property}: any family of pairwise intersecting closed balls has a non-empty global intersection.

\mk

Examples of geodesic injective metric spaces are normed vector spaces with the $\ell^\infty$ norm, and also finite-dimensional CAT(0) cube complexes with the piecewise $\ell^\infty$ metric (see~\cite{bowditch_median_injective}).

\mk

One key feature of the theory is that any metric space $X$ embeds isometrically in a unique minimal injective metric space, called the \emph{injective hull} of $X$ and denoted $EX$ (see~\cite{isbell}).

\mk

A metric space $(X,d)$ is called \emph{coarsely injective} if there exists a constant $C \geq 0$ such that, for any family $(x_i)_{i \in I}$ of points in $X$ and $(r_i)_{i \in \N}$ of nonnegative real numbers satisfying
$$ \forall i,j \in I, r_i+r_j \geq d(x_i,x_j),$$
the family of balls $(B(x_i,r_i+C))_{i \in \N}$ has a non-empty global intersection.

\mk

There is also a discrete version of injective metric spaces concerning graphs: a connected graph is called a \emph{Helly graph} if the family of combinatorial balls satisfy the Helly property: any family of pairwise intersecting balls has a non-empty global intersection.

\mk

Concerning actions of groups on injective metric spaces, we will distinguish three families:
\bit
\item A group $G$ is called \emph{coarsely injective} if it acts properly and coboundedly by isometries on an injective metric space, or equivalently it acts properly and cocompactly by isometries on a coarsely injective metric space (see~\cite[Proposition~3.12]{helly_groups}).
\item A group $G$ is called \emph{metrically injective} if it acts properly and cocompactly by isometries on an injective metric space.
\item A group $G$ is called \emph{Helly} if it acts properly and cocompactly by automorphisms on a Helly graph.
\eit

Any Helly group is metrically injective, by considering the injective hull of a Helly graph. And obvisouly, any metrically injective group is coarsely injective.

\mk

We now list examples of such groups.

\mk

According to~\cite{bandelt_vandevel_superextensions} (see also~\cite[Corollary~3.6]{hruskawise:packing}), the thickening of any CAT(0) cube complex is a Helly graph: in particular, any group acting properly and cocompactly on a CAT(0) cube complex is Helly. More generally, any group acting properly and cocompactly on a finite rank metric median space is metrically injective (see~\cite{bowditch_median_injective}). Urs Lang motivated the interest in group actions on injective metric spaces in~\cite{lang}, notably proving that any Gromov-hyperbolic group is Helly (see also~\cite{chepoi_estellon_packing}), and acts properly and cocompactly on the injective hull of any Cayley graph. Chalopin et al. proved (see~\cite[Corollary~6.2]{helly_groups}) that any type-preserving uniform lattice in a Euclidean building of type $\tilde{C}_n$ is Helly. Huang and Osajda proved that any Artin group of type FC is Helly (see~\cite{huang_osajda_helly}).

\mk

The author, Hoda and Petyt proved in~\cite{haettel_hoda_petyt} that any hierarchically hyperbolic group, including any mapping class group of a surface, is coarsely injective.

\mk

The existence of such actions on injective metric spaces enables us to deduce many properties reminiscent of non-positive curvature, let us list some of them:

\bthm
Assume that a finitely generated group $G$ is coarsely injective. Then:
\bit
\item $G$ is semi-hyperbolic in the sense of Alonso-Bridson, which has many consequences (\cite{bridson_haefliger}).
\item $G$ has finitely many conjugacy classes of finite subgroups (\cite[Proposition~1.2]{lang}).
\item $G$ satisfies the coarse Baum-Connes conjecture (\cite[Theorem~1.5]{helly_groups}).
\item Asymptotic cones of $G$ are contractible (\cite[Theorem~1.5]{helly_groups}).
\eit
Assume furthermore that $G$ is metrically injective. Then:
\bit
\item $G$ admits an EZ-boundary (\cite[Theorem~1.5]{helly_groups}).
\item $G$ satisfies the Farrell-Jones conjecture (see~\cite{kasprowski_rueping}).
\eit
Assume in addition that $G$ is a Helly group. Then:
\bit
\item $G$ is biautomatic (\cite[Theorem~1.5]{helly_groups}).
\eit
\ethm

Note that all consequences are already known for CAT(0) groups, except the biautomaticity (which does not hold for all CAT(0) groups, see~\cite{leary_minasyan}).

\mk

However, not all non-positively curved groups are coarsely injective: for instance, Hoda proved that the $(3,3,3)$ triangle Coxeter group, which is virtually $\Z^2$, is not Helly (see~\cite{hoda:crystallographic}).

\section{An injective distance on the space of all norms} \label{sec:injective_distance}

Let $\K$ denote a field (or a division algebra) with an absolute value $|\cdot| : \K \ra e^H \cup \{0\}$, where $H$ is a non-zero additive subgroup of $\R$. Let $V$ denote a $\K$-vector space. Recall that a norm on $V$ is a map $\eta : V \ra e^H$ that satisfies the following.
\bit
\item $\forall v \in V, \eta(v)=0 \Longleftrightarrow v=0$.
\item $\forall v \in V, \forall \alpha \in \K, \eta(\alpha v) = |\alpha| \eta(v)$.
\item $\forall u,v \in V, \eta(u+v) \leq \eta(u) + \eta(v)$.
\eit
Note that there is a natural partial order on the set of all norms on $V$: we say that $\eta \leq \eta'$ if $\forall v \in V, \eta(v) \leq \eta'(v)$. If $\eta \leq \eta'$, let us denote the interval $I(\eta,\eta')$ as the set of all norms $\theta$ such that $\eta \leq \theta \leq \eta'$.

\bpro \label{pro:helly_norms}
Let $X$ denote a non-empty set of norms on $V$ satisfying the following properties.
\bit
\item for every $\eta \in X$ and every $a \in H$, we have $e^a \eta \in X$.
\item for every $\eta,\eta' \in X$, there exist $a \in H$ such that $e^{-a} \eta' \leq \eta \leq e^a\eta'$.
\item the set $X$ is a join-semilattice: for every non-empty subset $F \subset X$ such that there exists $\eta \in X$ with $F \leq \eta$, the set $\{\eta' \in X \st F \leq \eta'\}$ has a unique minimum $\vee F \in X$.
\eit
For any two elements $\eta,\eta'$ in $X$, let us define the Goldman-Iwahori distance
$$d(\eta,\eta') = \sup_{v \in V \bs \{0\}} \left| \log \f{\eta(v)}{\eta'(v)} \right|.$$
Then the family of closed balls in the metric space $(X,d)$ satisfies the Helly property.
\epro

\bp
We will first describe balls in $(X,d)$. Fix $\eta \in X$ and $a \in \R_+$. Then $\eta' \in B(\eta,a)$ if and only if, for every $v \in V \bs \{0\}$, we have $-a \leq \log \f{\eta'(v)}{\eta(v)} \leq a$, hence $e^{-a} \eta(v) \leq \eta'(v) \leq e^a \eta(v)$. As a consequence, the ball $B(\eta,a)$ coincides with the interval $I(e^{-a} \eta, e^a \eta)$.

\mk

We will now prove that the intervals in $X$ satisfy the Helly property. Consider a family $(I_s=I(\eta_s,e^{2a_s} \eta_s))_{s \in S}$ of pairwise intersecting intervals in $X$, where $a_s \in H$ for each $s \in S$. Let $F=\{\eta_s\}_{s \in S} \subset X$: for any $s,t \in S$, since $I_s$ and $I_t$ are intersecting, we have $\eta_t \leq e^{2a_s} \eta_s$. According to the assumption on $X$, we can consider the join $\eta = \vee F \in X$. For each $s,t \in S$, since $\eta_t \leq e^{2a_s} \eta_s$, we deduce that $\eta \leq e^{2a_s} \eta_s$. In particular, for each $s \in S$, we have $\eta_s \leq \eta \leq e^{2a_s} \eta_s$, so $\eta \in I_s$. We have proved that the global intersection $\ds \bigcap_{s \in S} I_s$ is non-empty.
\ep

\section{Bruhat-Tits (extended) buildings are injective} \label{sec:buildings}

We will now apply Proposition~\ref{pro:helly_norms} to define an injective metric on classical Bruhat-Tits buildings.

\subsection{The standard and extended Bruhat-Tits buildings of $\GL(n,\K)$}

Let $\K$ be a field, with a non-Archimedean absolute value $|\cdot| : \K \ra \R_+$. Assume that $\K$ is a local field, or more generally that $\K$ is spherically complete: any decreasing intersection of balls in $\K$ has non-empty intersection. Let $V$ denote a $n$-dimensional vector space over $\K$.

\mk

Let us say that a map $\eta : V \ra \R_+$ is an ultrametric norm on $V$ if it satisfies the following.
\bit
\item $\forall v \in V, \eta(v)=0 \Longleftrightarrow v=0$.
\item $\forall v \in V, \forall \alpha \in \K, \eta(\alpha v) = |\alpha| \eta(v)$.
\item $\forall u,v \in V, \eta(u+v) \leq \max(\eta(u),\eta(v))$.
\eit

An ultrametric norm $\eta$ on $V$ is called \emph{diagonalizable} if there exists a basis $(v_1,\dots,v_n)$ of $V$ such that
$$\forall v=\sum_{i=1}^n x_i v_i \in V, \eta(v) = \max_{1 \leq i \leq n} |x_i|.$$

\mk

According to~\cite[Proposition~1.20]{remy_thuillier_werner}, if $\K$ is a local field, any ultrametric norm on $V$ is diagonalizable. This holds more generally if $\K$ is spherically complete, see~\cite[Remark~1.24]{remy_thuillier_werner}.

\mk

Say that two ultrametric norms $\eta,\eta' : V \ra \R_+$ are homothetic if there exists $a \in \R$ such that $\eta'=e^a \eta$. The set $\ov{X}$ of homothety classes of ultrametric norms on $V$ is called the Bruhat-Tits building of $\SL(n,\K)$ (see~\cite{parreau_immeubles} for instance).

\mk

Let $X$ denote the space of all (diagonalizable) ultrametric norms on $V$, it has been studied by Goldman and Iwahori (see~\cite{goldman_iwahori}) and can be identified with the extended Bruhat-Tits building of $\GL(n,\K)$. It is homeomorphic to the product $\ov{X} \times \R$.

\mk

For any two elements $\eta,\eta'$ in $X$, let us define the Goldman-Iwahori distance
$$d(\eta,\eta') = \sup_{v \in V \bs \{0\}} \left| \log \f{\eta(v)}{\eta'(v)} \right|.$$

We have an explicit description of the distance $d$ in terms of apartments of $X$. This description can also be found in~\cite{goldman_iwahori} without the building point of view, but we will give here a simple description using the building.

\mk

Let us recall the description of apartments in the Bruhat-Tits building $\ov{X}$ of $\GL(n,\K)$. For each basis $v_1,\dots,v_n$ of $V$ (up to homotheties and permutations), there is an associated apartment in $\ov{X}$. For each $m \in \R^n$, let us consider the following ultrametric norm on $V$:
$$\forall v=\sum_{i=1}^n x_i v_i \in V, \eta_m(v) = \max_{1 \leq i \leq n} e^{m_i} |x_i|.$$
Then the set of such homothety classes identifies with $\{x \in \R^n \st x_1+x_2+\dots+x_n=0\} \simeq \R^{n-1}$. It is a model of the standard Euclidean apartment of type $\tilde{A_{n-1}}$.

\mk

Let us now describe the apartments of the extended Bruhat-Tits building $X$ of $\GL(n,\K)$. For each basis $v_1,\dots,v_n$ of $V$ (up to homotheties and permutations), there is an associated apartment in $X$: the set of all norms $\{\eta_m \st m \in \R^n\}$ identifies with $\R^n$, which is a model of the extended Euclidean apartment of type $\tilde{A_{n-1}}$.

\bpro \label{pro:metric_infinity_apartment}
The metric $d$ on $X$ coincides with the $\ell^\infty$ metric on each extended apartment.
\epro

\bp
Let us denote by $d_\infty :X \times X \ra \R_+$ the map which to any couple $(x,y)$ in some apartment $A$ associates their $\ell^\infty$ distance in $A$. Note that $d_\infty$ is well-defined, but it is not obvious that it is a metric.

\mk

Fix a basis $v_1,\dots,v_n$ of $V$, and the associated apartment $A=\{\eta_m, m \in \R^n\}$ in $X$. Fix any $m \in \R^n$. Let $1 \leq i \leq n$ such that $|m_i| = \|m\|_\infty$, then we have 
$$\left|\log \f{\eta_m(v_i)}{\eta_0(v_i)}\right| = \left|\log e^{m_i}\right|=|m_i| = \|m\|_\infty,$$
hence $d_\infty(\eta_0,\eta_m) = \|m\|_\infty \leq d(\eta_0,\eta_m)$.

On the other hand, for any $v=\sum_{i=1}^n x_i v_i \in V$, we have 
\beq \left|\log \f{\eta_m(v)}{\eta_0(v)}\right| &=& \left|\log \f{\max_{1 \leq i \leq n} e^{m_i} |x_i|}{\max_{1 \leq i \leq n} |x_i|}\right| \\
&\leq& \left|\log \f{\max_{1 \leq i \leq n} e^{\|m\|_\infty} |x_i|}{max_{1 \leq i \leq n} |x_i|}\right| =\|m\|_\infty,\eeq
so we deduce that $d(\eta_0,\eta_m) \leq d_\infty(\eta_0,\eta_m)$.

So we have proved that $d(\eta_0,\eta_m) = d_\infty(\eta_0,\eta_m)$, for any $m \in \R^n$. Hence we deduce that $d=d_\infty$. 
\ep

We can now apply Proposition~\ref{pro:helly_norms} to prove that the metric $d$ is injective.

\bthm \label{thm:BTbuilding_injective}
The extended Bruhat-Tits building $X$ of $\GL(n,\K)$, endowed with the metric $d$, is injective.
\ethm

\bp
We first have to check that $X$ satisfies the three assumptions of Proposition~\ref{pro:helly_norms}.

\bit
\item For every $\eta \in X$ and every $a \in \R$, we know that $e^a \eta$ is an ultrametric norm on $V$, hence $e^a \eta \in X$.
\item For every $\eta,\eta' \in X$, let $a=d(\eta,\eta') = \sup_{v \in V \bs \{0\}} \left| \log \f{\eta(v)}{\eta'(v)} \right| \in \R_+$. For each $v \in V$, we have $\eta(v) \leq e^a \eta'(v)$ and $\eta'(v) \leq e^a \eta(v)$, hence $e^{-a} \eta' \leq \eta \leq e^a \eta'$.
\item For every non-empty subset $F \subset X$ such that there exists $\eta \in X$ with $F \leq \eta$, let $\theta = \sup F$. It it cleat that $\theta$ is a well-defined norm on $V$, we will check that it is ultrametric: fix $u,v \in V$. For every $\eps>0$, there exists $\eta' \in F$ such that $\theta(u+v) \leq \eta'(u+v)+\eps$. Then
$$\theta(u+v) \leq \eta'(u+v)+\eps \leq \max(\eta'(u),\eta'(v))+\eps \leq \max(\theta(u),\theta(v))+\eps.$$
This holds for any $\eps>0$, hence $\theta(u+v) \leq \max(\theta(u),\theta(v))$. So $\theta$ is an ultrametric norm on $V$: $\theta \in X$, and it is the unique minimum of the set $\{\eta' \in X \st F \leq \eta'\}$. Also recall that, since $\K$ is spherically complete, any ultrametric norm on $V$ is diagonalizable.
\eit

\mk

According to Proposition~\ref{pro:helly_norms}, the balls in $(X,d)$ satisfy the Helly property.

\mk

We also know by Proposition~\ref{pro:metric_infinity_apartment} that the metric space $(X,d)$ is geodesic. So we deduce that the metric space $(X,d)$ is injective.
\ep

\subsection{Case of a discrete valuation}

We will show that, if we further assume that the valuation is discrete, we can improve Theorem~\ref{thm:BTbuilding_injective} by finding a Helly graph.

\mk

Assume now that the absolute value is discrete: $|\cdot|(\K) = q^\Z \subset \R_+$, where $q$ is the cardinality of the residue field. Then the Bruhat-Tits building $\ov{X}$ of $\GL(n,\K)$ has a natural simplicial structure, where the vertex set $\ov{X}^{(0)}$ is given by the homothety classes of ultrametric norms with values in $q^\Z$.

\mk

Similarly, the extended Bruhat-Tits building $X$ of $\GL(n,\K)$ has a natural simplicial structure, where the vertex set $X^{(0)}$ is given by the ultrametric norms with values in $q^\Z$. To be consistent, we will in this case define the metric $d$ on $X^{(0)}$ as
$$d(\eta,\eta') = \sup_{v \in V \bs \{0\}} \left| \log_q \f{\eta(v)}{\eta'(v)} \right| \in \N.$$

\mk

Let us define the thickening $X'$ of $X$ as the graph with vertex set $X^{(0)}$, and with an edge between two vertices $\eta,\eta'$ if they satisfy $d(\eta,\eta')=1$.

\bthm \label{thm:gln_helly}
The thickening $X'$ of the extended Bruhat-Tits building of $\GL(n,\K)$ is a Helly graph.
\ethm

\bp
Following the same proof as Theorem~\ref{thm:BTbuilding_injective}, with $H=\log(q)\Z$, we prove that the integer-valued metric space $(X^{(0)},d)$ has the Helly property for balls.

\mk

It now suffices to prove that the distance $d$ is a graph distance. According to Proposition~\ref{pro:metric_infinity_apartment}, on each extended apartment, the metric $d$ coincides with the standard $\ell^\infty$ metric on $\R^n$. Since the restriction of the $\ell^\infty$ metric on $\R^n$ to the vertex set $\Z^n$ is a graph distance, we deduce that $d$ is a graph distance on $X^{(0)}$. This proves that the thickening $X'$ is a Helly graph.
\ep

\subsection{Classical Euclidean buildings}

We now show how to apply the previous results concerning the general linear group to the other classical groups.

\mk

Fix a local non-Archimedean field $\K$ with residual characteristic different from $2$, and consider a classical connected semisimple group $G$ over $\K$, realized as the identity component of the fixed point set of an involution $\Phi$ in a general linear group $\GL(n,\K)$. According to Bruhat and Tits (see~\cite{bruhat_tits_schemas} and \cite{prasad_yu}), the Bruhat-Tits building $X$ of $G$ identifies with the set of $\Phi$-fixed points in the Bruhat-Tits extended building $Y$ of $\GL(n,\K)$.

\mk

More generally, we may consider a finite group $F$ of automorphisms of $\GL(n,\K)$ such that the residual characteristic of $\K$ does not divide the order of $F$. Then, according to~\cite{prasad_yu}, the Bruhat-Tits building $X$ of $G=(\GL(n,\K)^F)^o$ identifies with the $F$-fixed points in the Bruhat-Tits extended building $Y$ of $\GL(n,\K)$.

\mk

Endow $X$ with the induced piecewise $\ell^\infty$ metric $d$ from $Y$.

\bthm \label{thm:classical_injective}
The Bruhat-Tits building $X$ of $G$, with the metric $d$, is injective.
\ethm

\bp
According to~\cite[Proposition~1.2]{lang}, the fixed point set $X=Y^F$ of any finite group action on an injective metric space is non-empty and injective. So the metric space $(X,d)$ is injective. 
\ep

We can also strengthen this result by looking for an action of $G$ on a Helly graph.

\bthm \label{thm:classical_helly}
The group $G$ acts properly and cocompactly by automorphisms on a Helly graph.
\ethm

\bp
Let $Y'$ denote the thickening of the $0$-skeleton of $Y$, which is a Helly graph according to Theorem~\ref{thm:gln_helly}. Let $F(Y')$ denote the face complex of $Y'$: it is the simplicial complex with vertex set the set of cliques of $Y'$, and with simplices the set of cliques contained in a given clique of $Y'$. According to~\cite[Lemma~5.30]{helly_groups}, the face complex $F(Y')$ is clique-Helly (i.e. the family of maximal cliques satisfies the Helly property).

\mk

The group $\GL(n,\K)$ acts properly and cocompactly on $Y'$. Let $X'$ denote the fixed point set of $F$ inside $F(Y')$: according to~\cite[Theorem~7.1, Corollary~7.4]{helly_groups}, it is a non-empty clique-Helly graph. According to~\cite{chalopin_chepoi_hirai_osajda}, the underlying graph of $X'$ is Helly, and $G$ acts properly and cocompactly on $X'$.
\ep

The following is immediate.

\bcor \label{cor:classical_padic}
Let $G$ denote a classical reductive Lie group over a non-Archimedean local field of characteristic different from $2$, and let $a \geq 0$ denote the number of semisimple factors of type $A$. Then $G \times \Z^a$ acts properly and cocompactly by automorphisms on a Helly graph.

For any cocompact lattice $\Gamma$ in $G$, the group $\Gamma \times \Z^a$ acts properly and cocompactly by automorphisms on a Helly graph, and the group $\Gamma$ is biautomatic.\ecor

\bp
This is a direct consequence of Theorem~\ref{thm:classical_helly}. According to~\cite[Theorem~1.5]{helly_groups}, any Helly group is biautomatic. And according to~\cite[Theorem~B]{mosher_biautomatic}, every direct factor of a biautomatic group is biautomatic.
\ep

Swiatkowski proved that any group acting properly and cocompactly on any Euclidean building is biautomatic (see~\cite[Theorem~6.1]{swiatkowski_biautomatic}). So we obtain another point of view on this result, for uniform lattices in classical groups.

\section{Symmetric spaces are coarsely injective} \label{sec:symmetric_spaces}

We will use Proposition~\ref{pro:helly_norms} to find the injective hull of the symmetric space of $\GL(n,\R)$, and to study the injective hulls of classical symmetric spaces of non-compact type.

\subsection{The symmetric space of $\GL(n,\R)$}

Fix $\K=\R, \C$ or $\H$ (the division algebra of quaternions), fix $n \geq 2$, and let $V$ denote a $n$-dimensional vector space over $\K$.

Say that two Euclidean norms $\eta,\eta' : V \ra \R_+$ are homothetic if there exists $a \in \R$ such that $\eta'=e^a \eta$. The set $\ov{X}$ of homothety classes of hermitian norms on $V$ is called the symmetric space of $\SL(n,\K)$, and it identifies naturally with the homogeneous space $\SL(n,\K) / \SU(n,\K)$.

\mk

Let $X$ denote the space of all hermitian norms on $V$, it is called the symmetric space of $\GL(n,\K)$ and it identifies naturally with the homogeneous space $\GL(n,\K) / \U(n,\K)$. It is homeomorphic to the product $\ov{X} \times \R$.

\mk

Let $\hat{X}$ denote the space of all norms on $V$ which are invariant under the unit group $\mathbb{U}$ of $\K$, it contains $X$ as the subset of hermitian norms. The space $\hat{X}$ can also be described as the space of all compact convex subsets of $V$ with non-empty interior, which are invariant under the linear diagonal action of the unit group $\mathbb{U}$ of $\K$. Such convex subsets will be called \emph{symmetric}. We will call it the augmented symmetric space of $\GL(n,\K)$. The group $\GL(n,\K)$ acts naturally on $\hat{X}$, by precomposing the norms, or by the linear action on convex subsets of $V$.

\mk

For any two elements $\eta,\eta'$ in $\hat{X}$, let us define the distance
$$\hat{d}(\eta,\eta') = \sup_{v \in V \bs \{0\}} \left| \log \f{\eta(v)}{\eta'(v)} \right|.$$
It is a lift of the Banach-Mazur distance, which is defined on the set of isometry classes of such norms.

\mk

Let us also define the distance $d$ on $X$ as the restriction of the distance $\hat{d}$.

\mk

We have an explicit description of the distance $d$ in terms of maximal flats of $X$.

\mk

Let us recall the description of maximal flats in the symmetric space $\ov{X}$ of $\SL(n,\K)$. For each basis $v_1,\dots,v_n$ of $V$ (up to homotheties and permutations), there is an associated maximal flat in $\ov{X}$. For each $m \in \R^n$, let us consider the following hermitian norm on $V$:
$$\forall x=\sum_{i=1}^n x_i v_i \in V, \eta_m(x) = \sqrt{\sum_{i=1}^n e^{2m_i} |x_i|^2}.$$
Then the set of such homothety classes identifies with $\{m \in \R^n \st m_1+m_2+\dots+m_n=0\} \simeq \R^{n-1}$. It is a model of the standard Euclidean flat of type $\tilde{A_{n-1}}$.

\mk

Let us now describe the maximal flats of the symmetric space $X$ of $\GL(n,\K)$. For each basis $v_1,\dots,v_n$ of $V$ (up to homotheties and permutations), there is an associated maximal flat in $X$, the set $\{\eta_m \st m \in \R^n\}$ is a model of the extended Euclidean flat of type $\tilde{A_{n-1}}$.

\bpro \label{pro:real_geodesic}
The metric $d$ on $X$ coincides with the $\ell^\infty$ metric on each maximal flat.\epro

\bp
Let us denote by $d_\infty :X \times X \ra \R_+$ the map which to any couple $(x,y)$ in some maximal flat $A$ associates their $\ell^\infty$ distance in $A$. As in Proposition~\ref{pro:metric_infinity_apartment}, we prove that $d=d_\infty$.
\ep

\bpro \label{pro:real_bounded}
The symmetric space $X$ of $\GL(n,\K)$ is cobounded in $\hat{X}$.
\epro

\bp
Let $K \in \hat{X}$. Let $B \subset K$ denote the unique John-L\"owner ellipsoid of maximal volume. Since $K$ is invariant under the linear diagonal action of the unit group $\mathbb{U}$, by uniqueness of $B$, we deduce that $B$ is also invariant under the linear diagonal action of the unit group $\mathbb{U}$. So the convex $B$ is the unit ball of a hermitian norm on $\K^n$: $B \in X$. According to~\cite{john}, we know that $\hat{d}(B,K) \leq \log(\sqrt{an})$, where $a=\dim_\R(\K)$. Therefore any point of $\hat{X}$ is at distance at most $\log(\sqrt{an})$ from $X$.
\ep

We could then apply directly Proposition~\ref{pro:helly_norms} to deduce that balls in $\hat{X}$ satisfy the Helly property. However, it is not clear yet that $\hat{X}$ is a geodesic metric space. Moreover, it is interesting to describe explicitly the injective hull of $X$. So instead of using Proposition~\ref{pro:helly_norms}, we will prove directly that $\hat{X}$ is the injective hull of $X$.

\bthm \label{thm:real_coarse_helly}
Let $X$ denote the symmetric space of $\GL(n,\K)$, endowed with the $\ell^\infty$ distance $d$. The injective hull of $X$ is the space $\hat{X}$ of all symmetric compact convex subspaces of $\K^n$ with non-empty interior. Moreover, $X$ is cobounded in $\hat{X}$, which is is proper. As a consequence, $\GL(n,\K)$ acts properly and cocompactly on the injective space $\hat{X}$.
\ethm

\bp
We will use Lang's description of the injective hull of $(X,d)$ (see~\cite{lang}). Let us denote
$$\Delta X=\{f : X \ra \R_+ \mbox{ $1$-Lipschitz} \st \forall x,y \in X, f(x)+f(y) \geq d(x,y)\},$$
equipped with the supremum metric:
$$\forall f,f' \in \Delta X, d_{\Delta X}(f,f')=\sup_{x \in X} |f(x)-f'(x)|.$$
Let us denote by $EX$ the set of minimal elements of $(\Delta X,\leq)$. More explicitely, we have
$$EX = \{f \in \Delta X \st \forall x \in X, f(x) = \sup_{y \in X} d(x,y)-f(y)\}.$$
There is a canonical isometric embedding $e:X \ra EX$ defined by $x \mapsto d(x,\cdot)$, and $EX$ is the injective hull of $X$.

\mk

We will now define an isometric embedding $\phi$ from $\hat{X}$ into $EX$ extending $e$. For each convex subset $C \in \hat{X}$, let us consider $\phi(C) : B \in X \mapsto \hat{d}(C,B)$: it is clear that $\phi(C) \in \Delta X$.

\mk

We will prove that, for any $C,C' \in \hat{X}$ and any $\eps>0$, there exists $B \in X$ such that $\hat{d}(B,C')+\hat{d}(C',C) \leq \hat{d}(B,C)+\eps$. Fix $C,C' \in \hat{X}$ distinct, and let $t=\hat{d}(C,C')>0$. Without loss of generality, we may assume that for every $s<t$, we have $C \not\subset e^sC'$. Let $v \in \partial e^{-t}C \cap \partial C'$.

Fix $\eps>0$. There exists an ellipsoid $B' \in X$ such that $C' \subset B'$ and $v \in \partial e^{-\eps} B'$. Fix $a>0$ large enough such that $B=e^{-a}B' \subset C \cap C'$. Then $\hat{d}(B,C') \leq a$, and since $v \in \partial e^{a-\eps}B \cap \partial e^{-t}C$, we deduce that $\hat{d}(B,C) \geq a+t-\eps$. Hence we have $\hat{d}(B,C')+\hat{d}(C',C) \leq a+t \leq \hat{d}(B,C)+\eps$.

\mk

In particular, this result implies that the map $C \in \hat{X} \mapsto \phi(C) \in \Delta X$ is an isometric embedding. Furthermore, for any $C \in \hat{X}$ and $C' \in X$, according to the same result, we deduce that $\phi(C)(C')=\sup_{B \in X} \hat{d}(B,C')-\phi(C)(B)$, hence $\phi(C) \in EX$.

\mk

We will now prove that $\phi$ extends $e$: for any $B,C \in X$, we have $\phi(C)(B)=\hat{d}(C,B)=d(C,B)=e(C)(B)$. So we have proved that $\phi$ is an isometric embedding of $\hat{X}$ into $EX$, extending $e$.

\mk

To conclude, we will prove that $\phi$ is surjective: let $f \in EX$, and consider $C = \bigcap_{B \in X} e^{f(B)}B$. Fix $B_0 \in X$, then $C \subset e^{f(B_0)}$. On the other hand, for any $B \in X$, we have $e^{-f(B_0)}B_0 \subset e^{f(B)}B$, hence $e^{-f(B_0)}B_0 \subset C$. We deduce that $\hat{d}(C,B_0) \leq f(B_0)$, for any $B_0 \in X$. Hence $\phi(C) \leq f$, and by minimality of $f \in EX$ we conclude that $\phi(C)=f$. So $\phi$ is surjective.

\mk

We have proved that $(\hat{X},\hat{d})$ is isometric to $EX$, hence it is the injective hull of $X$.

\mk

According to Proposition~\ref{pro:real_bounded} we know that $X$ is cobounded in $\hat{X}$, which is locally compact. Since $\hat{X}$ is also complete and geodesic, we deduce that $\hat{X}$ is proper.
\ep

\subsection{Classical symmetric spaces of non-compact type}

We now show how to apply the previous results concerning the general linear group to the other classical groups.

\mk

Say that a semisimple non-compact real Lie group $G$ over $\R$ is \emph{classical not of type $\SL$} if it is commensurable to one of $\Sp(n,\R)$, $\Sp(n,\C)$, $\Sp(n,\H)=\SU^*(2n)$, $\O(n,\C)$, $\O(n,\H)=\SO^*(2n)$, $\O(p,q)$, $\U(p,q)$, $\Sp(p,q)$ (see~\cite{helgason}). Note that, due to some exceptional isomorphisms (such as $\SL(4,\R)$ being commensurable to $\SO(3,3)$), such a group $G$ may also be commensurable to a group $\SL(n,\K)$.

\mk

There exists $n \geq 1$, $\K=\R,\C$ or $\H$, and a finite group $F$ of automorphisms of $\GL(n,\K)$ such that $G$ embeds in $\GL(n,\K)$ and identifies with the fixed point subgroup $\GL(n,\K)^F$. More explicitely, $F$ is generated by the involution $A \mapsto J^{-1}(A^*)^{-1}J$, where $J$ is the matrix associated with the form defining $G$.

Furthermore, if we denote by $K$ a maximal compact subgroup of $G$, we can assume that $K=\U(n)^F$, and that the corresponding embedding of the symmetric space $X=G/K$ of $G$ into the symmetric space $Y=\GL(n,\K)/\U(n)$ has image the fixed point set $X=Y^F$ of $F$. We endow $X$ with the restriction of the $\ell^\infty$ length metric on $Y$. Let us denote $\hat{Y}$ the space of all symmetric compact convex subspaces of $\K^n$ with non-empty interior. According to Theorem~\ref{thm:real_coarse_helly}, $\hat{Y}$ is also the injective hull of $Y$. Let us denote $\hat{X}=\hat{Y}^F$.

\bpro
Any classical irreducible symmetric space of non-compact type $X$, which is not of type $\SL$, is cobounded in $\hat{X}$.
\epro

\bp
Let $K \in \hat{Y}^F$. Let $B \subset K$ denote the unique John-L\"owner ellipsoid of maximal volume. By uniqueness, we deduce that $B$ is invariant under $F$, and also under the unit group $\U$ of $\K$, i.e. $B \in X=Y^F$.  According to~\cite{john}, we know that $d(B,K) \leq \log(\sqrt{an})$, where $a=\dim_\R(\K)$. Therefore any point of $\hat{X}=\hat{Y}^F$ is at distance at most $\log(\sqrt{an})$ from $X=Y^F$.
\ep

\bthm \label{thm:classical_coarse_helly}
Let $X$ denote a classical irreducible symmetric space of non-compact type which is not of type $\SL$. Then the Finsler metric space $(X,d)$ is coarsely injective, and its injective hull is proper.
\ethm

\bp
According to Theorem~\ref{thm:real_coarse_helly}, the symmetric space $Y=\GL(n,\K)/\U(n)$, endowed with the piecewise $\ell^\infty$ distance, is coarsely injective, and its injective hull $\hat{Y}$ is proper. The isometric action of the finite group $F$ on $Y$ extends to an isometric action on $\hat{Y}$.

\mk

According to~\cite[Proposition~1.2]{lang}, the fixed point set $\hat{Y}^F$ of $F$ on $\hat{Y}$ is an injective metric space. Therefore, the injective hull $EX$ of $X$ may be realized as an isometric closed subspace of $\hat{X}=\hat{Y}^F$, so $EX$ is proper.

\mk

On the other hand, since $X=Y^F$ is cobounded in $\hat{X}=\hat{Y}^F$, we deduce that $X$ is cobounded in $EX$.
\ep

Note that if $X$ has rank $1$, we have a similar result.

\bpro
Let $X$ denote a rank $1$ symmetric space of non-compact type, and let $d$ denote the standard Riemannian metric on $X$. Then the metric space $(X,d)$ is coarsely injective, and its injective hull is proper.
\epro

\bp
The metric $d$ is Gromov-hyperbolic, so according to~\cite[Proposition~1.3]{lang} we know that $(X,d)$ is coarsely injective. We also know that $X$ may be realized as a totally geodesic subspace of the symmetric space $Y$ of $\GL(n,\R)$ (for the Riemannian metric on $Y$) for some $n \geq 2$ (see~\cite[Theorem~1.6.5]{eberlein}). Let us denote by $d_Y$ the $\ell^\infty$ metric on $Y$. Since each Riemannian geodesic is a $d_Y$ geodesic, we deduce that $d$ coincides with the restriction on $X$ of $d_Y$ (up to a constant factor, which may be chosen to be $1$). Hence $(X,d)$ is isometrically embedded in the proper injective metric space $E(Y,d_Y)$, so the injective hull of $X$ is proper.
\ep

The following consequence of Theorem~\ref{thm:classical_coarse_helly} is immediate.

\bcor \label{cor:classical_real_reductive}
Let $G$ denote any reductive Lie group over $\R$, with classical non-compact semisimple factors. Let $a \geq 0$ denote the number of semisimple factors of type $\SL$. Then $G \times \R^a$ acts properly and cocompactly on an injective metric space. In particular, for any cocompact lattice $\Gamma$ in $G$, the group $\Gamma \times \Z^a$ acts properly and cocompactly on an injective metric space.
\ecor

As we will see below, the factors $\R^a$ and $\Z^a$ are necessary.

\section{The special linear group is not coarsely injective} \label{sec:sln_not_helly}

We now turn to the case of the special linear group. We will prove that it is not coarsely injective, inspired by the result of Hoda that the $(3,3,3)$ triangle Coxeter group $W$, which is virtually $\Z^2$, is not Helly (see~\cite{hoda:crystallographic}). However, the group $W$ is a subgroup of $\Z^3 \rtimes \frak{S}_3$, which is Helly. This situation is analogous to the inclusion of $\SL(n,\K)$ in $\GL(n,\K)$:

\bthm \label{thm:sln_not_helly}
Let $\K$ be a local field (with characteric different from $2$ if $\K$ in non-Archimedean), and let $n=3$ or $n \geq 5$. Then $\SL(n,\K)$ is not coarsely injective: $\SL(n,\K)$ does not act properly and coboundedly on an injective metric space.
\ethm

Note that $\SL(4,\R)$ is commensurable to $\SO(3,3)$, so according to Theorem~\ref{thm:classical_coarse_helly} it is coarsely injective. We do not know about $\SL(4,\K)$, when $K \neq \R$.

\bp
By contradiction, assume that $G=\SL(n,\K)$ acts properly and coboundedly on an injective metric space $X$.

\mk

Let $A \subset \SL(n,\K)$ denote the diagonal subgroup, and let $M \subset \SL(n,\K)$ denote the monomial subgroup of $\SL(n,\K)$: $M \simeq A \rtimes \frak{A}_n$ is the subgroup of matrices with exactly one non-zero entry on each row and each column (and $\frak{A}_n$ denotes the alternating group). Let $F \subset A$ denote the finite diagonal subgroup with entries in $\{-1,1\}$. Since $\K$ has characteristic different from $2$, we know that the subgroup of $G$ fixed by the conjugation by $F$ is $G^F=A$. According to~\cite[Proposition~1.2]{lang}, the fixed point set $X^F$ of $F$ in $X$ is non-empty and injective. We will prove that $M$ acts properly and coboundedly on the injective metric space $X^F$. Firstly, since $F$ is normalized by $M$, we deduce that $M$ stabilizes $X^F$, and acts properly on $X^F$. We will prove that $A$ acts coboundedly on $X^F$, which will imply that $M$ also acts coboundedly on $X^F$.

\mk

Fix $x_0 \in X^F$, and let $C_X \geq 0$ such that any $x \in X$ is at distance at most $C_X$ from a point in $G \cdot x_0$.

\mk

Fix $x \in X^F$, there exists $g \in G$ such that $d(x,g \cdot x_0) \leq C_X$. So we deduce that, for any $f \in F$, we have $d(g \cdot x_0, fg \cdot x_0) \leq 2C_X$. Let $d_G$ denote a proper left-invariant metric on $G$. Since the action of $G$ on $X$ is proper, we deduce that there exists $C_G \geq 0$ such that, for any $f \in F$, we have $d_G(g,fgf^{-1}) \leq C_G$. Let $Y$ denote the symmetric space or Bruhat-Tits building of $G$, endowed with the CAT(0) metric, choose a basepoint $y_0 \in Y$ fixed by $F$, and let $y=g \cdot y_0$. Then there exists $C_Y \geq 0$ such that, for any $f \in F$, we have $d(y,f \cdot y) \leq C_Y$. Let $\ov{y} \in Y$ denote the CAT(0) barycenter of the finite orbit $F \cdot y$: it is fixed by $F$, and also $d(y,\ov{y}) \leq C_Y$. Since $G$ acts coboundedly on $Y$, there exists a constant $C'_G \geq 0$ and $\ov{g} \in G^F=A$ such that $d_G(g,\ov{g}) \leq C'_G$. Let us denote $\ov{x} = \ov{g} \cdot x_0 \in X^F$: there exists a constant $C'_X$ such that $d(\ov{x},x) \leq C'_X$. This proves that the action of $A$ on $X^F$ is cobounded.

\mk

So we have proved that the group $M \simeq A \rtimes \frak{A}_n$ acts properly and coboundedly on the injective metric space $X^F$. We deduce the existence of a proper, left-invariant metric $d_M$ on $M$ which is coarsely injective.

\mk

Let us consider a non-principal ultrafilter $\omega$ on $\N$, and a sequence $(\lambda_k)_{k \in \N}$ in $(0,\infty)$ which $\omega$-converges to $0$. Note that $A$ is isomorphic to $(\K^*)^{n-1} \simeq \{x \in (\K^*)^n \st x_1\times \dots  \times x_n=1\}$, where $\K^*$ denotes the multiplicative group of $\K$. Also note that any asymptotic cone of $\K^*$ is isomorphic to $(\R_+^*,\times)$, which we will realize as $(\R,+)$. Consider the asymptotic cone $(A_\infty,e_\infty,d_\infty) = \liml_{k \in \omega} (A,1,\lambda_k d_M)$, it is isomorphic to the group $\R^{n-1} \simeq \{x \in \R^n \st x_1+ \dots +x_n=0\}$. Since $A$ is abelian, we deduce that $A_\infty$ acts on itself by left translations: if $[a_k]_{k \in \N}=[a'_k]_{k \in \N} \in A_\infty$ and $[b_k]_{k \in \N}=[b'_k]_{k \in \N} \in A_\infty$, then
$$\liml_{k \in \omega} \lambda_kd_M(a_kb_k,a'_kb'_k) \leq \liml_{k \in \omega} \lambda_kd_M(a_k,a'_k)+\lambda_kd_M(b_k,b'_k) = 0.$$
This action of $A_\infty$ on itself preserves the metric metric $d_\infty$. Hence we deduce that $d_\infty$ is a norm on $A_\infty \simeq \R^{n-1}$. Also note that the natural action of $\frak{A}_n$ on $A$ induces the natural action of $\frak{A}_n$ on $A_\infty \simeq \R^{n-1} \simeq \{x \in \R^n \st x_1+ \dots +x_n=0\}$, and it is isometric with respect to $d_\infty$.

\mk

According to~\cite{nachbin}, the only $(n-1)$-dimensional injective normed vector spaces are linearly isometric to $\ell_\infty^{n-1}$. The linear isometry group of $\ell_\infty^{n-1}$ is the isometry group of the $(n-1)$-cube, $\frak{S}_{n-1} \ltimes \{\pm 1\}^{n-1}$. So we deduce that there exists an injective group morphism from $\frak{A}_n$ to $\frak{S}_{n-1} \ltimes \{\pm 1\}^{n-1}$. If $n=3$, then $\frak{A}_3$ has an order $3$ element and $\frak{S}_{2} \ltimes \{\pm 1\}^{2}$ does not, which is a contradiction. If $n \geq 5$, then $\frak{A}_n$ is simple, and there is no injective morphism from $\frak{A}_n$ to either $\frak{S}_{n-1}$ or $\{\pm 1\}^{n-1}$, which is a contradiction.

\mk

This concludes the proof that $\SL(n,\K)$ is not coarsely injective. 
\ep

However, this leaves the following question open: are uniform lattices in $\SL(n,\K)$ coarsely injective~?

\sign

\bibliographystyle{alpha}
\bibliography{../../../bibli}

\newcommand{\etalchar}[1]{$^{#1}$}
\def\polhk#1{\setbox0=\hbox{#1}{\ooalign{\hidewidth
  \lower1.5ex\hbox{`}\hidewidth\crcr\unhbox0}}}
\begin{thebibliography}{CCG{\etalchar{+}}20}

\bibitem[BH99]{bridson_haefliger}
Martin~R. Bridson and Andr\'e Haefliger.
\newblock {\em {Metric \,spaces \,of \,non-positive \,curva\-ture}}, volume
  {319\!} of {\em {Grund.~math.~Wiss.}}
\newblock {Springer}, 1999.

\bibitem[Bow20]{bowditch_median_injective}
Brian~H. Bowditch.
\newblock Median and injective metric spaces.
\newblock {\em Math. Proc. Cambridge Philos. Soc.}, 168(1):43--55, 2020.

\bibitem[BT84]{bruhat_tits_schemas}
F.~Bruhat and J.~Tits.
\newblock Sch\'{e}mas en groupes et immeubles des groupes classiques sur un
  corps local.
\newblock {\em Bull. Soc. Math. France}, 112(2):259--301, 1984.

\bibitem[BvdV91]{bandelt_vandevel_superextensions}
H.-J. Bandelt and M.~van~de Vel.
\newblock Superextensions and the depth of median graphs.
\newblock {\em J. Combin. Theory Ser. A}, 57(2):187--202, 1991.

\bibitem[CCG{\etalchar{+}}20]{helly_groups}
J{\'e}r{\'e}mie Chalopin, Victor Chepoi, Anthony Genevois, Hiroshi Hirai, and
  Damian Osajda.
\newblock Helly groups.
\newblock {\em arXiv preprint arXiv:2002.06895}, 2020.

\bibitem[CCHO21]{chalopin_chepoi_hirai_osajda}
J{\'e}r{\'e}mie Chalopin, Victor Chepoi, Hiroshi Hirai, and Damian Osajda.
\newblock Weakly modular graphs and nonpositive curvature.
\newblock {\em Mem. Amer. Math. Soc.}, 2021.

\bibitem[CE07]{chepoi_estellon_packing}
Victor Chepoi and Bertrand Estellon.
\newblock Packing and covering $\delta$-hyperbolic spaces by balls.
\newblock In {\em Approximation, {R}andomization, and {C}ombinatorial
  {O}ptimization. {A}lgorithms and {T}echniques}, pages 59--73. Springer, 2007.

\bibitem[Ebe96]{eberlein}
Patrick~B. Eberlein.
\newblock {\em {Geometry of nonpositively curved manifolds}}.
\newblock {Chicago Lect.~Math.} {The University of Chicago Press}, 1996.

\bibitem[GI63]{goldman_iwahori}
O.~Goldman and N.~Iwahori.
\newblock The space of p-adic norms.
\newblock {\em Acta Math.}, 109:137--177, 1963.

\bibitem[Hel78]{helgason}
Sigurdur Helgason.
\newblock {\em {Differential geometry, Lie groups, and symmetric spaces\!}}
\newblock {Grad.~Stud.~Math.~{\bf 34}, Amer.~Math.~Soc.} 1978.

\bibitem[HHP21]{haettel_hoda_petyt}
Thomas Haettel, Nima Hoda, and Harry Petyt.
\newblock {Coarse injectivity, hierarchical hyperbolicity, and
  semihyperbolicity}.
\newblock {\em to appear in Geom. Topol.}, 2021.

\bibitem[Hir20]{hirai_uniform_modular}
Hiroshi Hirai.
\newblock Uniform modular lattices and affine buildings.
\newblock {\em Adv. Geom.}, 20(3):375--390, 2020.

\bibitem[HO21]{huang_osajda_helly}
Jingyin Huang and Damian Osajda.
\newblock {Helly meets Garside and Artin}.
\newblock {\em Invent. Math.}, 225(2):395--426, 2021.

\bibitem[Hod20]{hoda:crystallographic}
Nima Hoda.
\newblock Crystallographic {H}elly groups.
\newblock {\em arXiv preprint arXiv:2010.07407}, 2020.

\bibitem[HW09]{hruskawise:packing}
G.~Christopher Hruska and Daniel~T. Wise.
\newblock Packing subgroups in relatively hyperbolic groups.
\newblock {\em Geom. Topol.}, 13(4):1945--1988, 2009.

\bibitem[Isb64]{isbell}
J.~R. Isbell.
\newblock Six theorems about injective metric spaces.
\newblock {\em Comment. Math. Helv.}, 39:65--76, 1964.

\bibitem[Joh48]{john}
Fritz John.
\newblock Extremum problems with inequalities as subsidiary conditions.
\newblock In {\em Studies and {E}ssays {P}resented to {R}. {C}ourant on his
  60th {B}irthday, {J}anuary 8, 1948}, pages 187--204. Interscience Publishers,
  Inc., New York, N. Y., 1948.

\bibitem[KR17]{kasprowski_rueping}
Daniel Kasprowski and Henrik R\"uping.
\newblock The {F}arrell-{J}ones conjecture for hyperbolic and {CAT}(0)-groups
  revisited.
\newblock {\em J. Topol. Anal.}, 9(4):551--569, 2017.

\bibitem[Lan13]{lang}
Urs Lang.
\newblock Injective hulls of certain discrete metric spaces and groups.
\newblock {\em J. Topol. Anal.}, 5(3):297--331, 2013.

\bibitem[LM21]{leary_minasyan}
Ian Leary and Ashot Minasyan.
\newblock {Commensurating HNN-extensions: non-positive curvature and
  biautomaticity}.
\newblock {\em Geom. Topol.}, 2021.

\bibitem[Mos97]{mosher_biautomatic}
Lee Mosher.
\newblock Central quotients of biautomatic groups.
\newblock {\em Comment. Math. Helv.}, 72(1):16--29, 1997.

\bibitem[Nac50]{nachbin}
Leopoldo Nachbin.
\newblock A theorem of the {H}ahn-{B}anach type for linear transformations.
\newblock {\em Trans. Amer. Math. Soc.}, 68:28--46, 1950.

\bibitem[Par99]{parreau_immeubles}
Anne Parreau.
\newblock {\em {Immeubles euclidiens : construction par les normes et \'etudes
  des isom\'etries}}.
\newblock \dans \og Cristallographic Groups and Their Generalisations, II \fg.
  Contemp. Math., Kortrijk, Amer. Math. Soc., 1999.

\bibitem[PY02]{prasad_yu}
Gopal Prasad and Jiu-Kang Yu.
\newblock On finite group actions on reductive groups and buildings.
\newblock {\em Invent. Math.}, 147(3):545--560, 2002.

\bibitem[RTW12]{remy_thuillier_werner}
Bertrand R\'emy, Amaury Thuillier, and Annette Werner.
\newblock {Bruhat-Tits Theory from Berkovich's Point of View. II - Satake
  compactifications.}
\newblock {\em {J.I.M.J.}}, 2(11):421--465, 2012.

\bibitem[\'{S}06]{swiatkowski_biautomatic}
Jacek \'{S}wi\.{a}tkowski.
\newblock Regular path systems and (bi)automatic groups.
\newblock {\em Geom. Dedicata}, 118:23--48, 2006.

\end{thebibliography}

\end{document}